\documentclass[12pt]{article}
\newtheorem{theorem}{Theorem}
\newtheorem{lemma}{Lemma}

\newtheorem{conjecture}{Conjecture}
\newtheorem{question}{Question}

\input amssym.tex
\newtheorem{preproof}{{\bf Proof.}}

\newenvironment{proof}[1]{\begin{preproof}{\rm
               #1}\hfill{\rule[-0.5mm]{2mm}{2mm}}}{\end{preproof}}
\def\n#1{\vbox to 3mm{\vspace{1mm}\vfill \hbox to 2.0mm{\hfill
             $#1$\hfill} \vfill }}
\def\arraystretch{1.0}                 
\def\lcs#1{{\rm lcs}(#1)}
\def\scs#1{{\rm scs}(#1)}

\usepackage{array}
\newlength{\cellwid}
\makeatletter
\newenvironment{latinsq}[1][00]{%
 \def\centcol{\centering \let\\=\tabularnewline \CellStrut}
 \arraycolsep0pt
 \setbox\@tempboxa\hbox{#1}%
 \cellwid\ht\@tempboxa  \advance\cellwid\dp\strutbox
 \ifdim\cellwid<\wd\@tempboxa 
   \@tempdima .5\wd\@tempboxa \advance\@tempdima -.5\cellwid 
   \advance\cellwid \@tempdima \advance\@tempdima\dp\strutbox
 \else
   \@tempdima\dp\strutbox
 \fi
 \edef\CellStrut{\vrule
   width\z@ height\the\cellwid depth\the\@tempdima \relax}
 \advance\cellwid\@tempdima \advance\cellwid-2\arraycolsep
 \array{|>{\centcol}p{\cellwid}|*{20}{>{\centcol}p{\cellwid}|}}}%
 {\endarray}
\makeatother

\title{{\bf A new bound on the size of the largest critical set in a Latin square}}
\author{\sc Richard Bean and E.S. Mahmoodian}
\date{}
\begin{document}
\maketitle
\begin{center}
Centre for Discrete Mathematics and Computing \\
Department of Mathematics \\
The University of Queensland \\
Queensland 4072, Australia \\
and\\
 Department of Mathematical Sciences \\
 Sharif University of Technology \\
 P.O. Box 11365--9415 \\
 Tehran, I.R. Iran
\end{center}
\begin{abstract}
A {\sf critical set }
in an $n \times n$ array is a set $C$ of given entries, such that
there exists a unique extension of  $C$ to  an
$n\times n$ Latin square and no proper subset of $C$ has this property.
The cardinality of the largest critical set in any Latin square
of order $n$ is denoted by $\lcs{n}$.  In 1978 Curran and van Rees proved that
$\lcs{n} \leq n^2 - n$.  Here we show that $\lcs{n} \leq n^2-3n+3$.
\end{abstract}
\section{Introduction}        
For the purposes of this paper, a {\sf Latin square} of order $n$ is an
$n$ $\times$ $n$ array of integers chosen from the set
$X = \{1,2, \ldots, n\}$ such that
each integer occurs exactly once in each row and exactly once in each column.
An example of a Latin square of order 4 is shown below.
$$
\begin{array}{|c|c|c|c|}
\hline 1&2&3&4 \\
\hline 2&3&4&1 \\
\hline 3&4&1&2 \\
\hline 4&1&2&3 \\
\hline \end{array}
$$
A Latin square can also be written as a set of ordered triples
$\{ (i,j;k) \mid$  symbol $k$ occurs in position $(i,j)$ of the array$\}$.

A {\sf partial Latin square} $P$ of order $n$ is an $n\times n$ array
with entries chosen
from the set $X = \{1,2, \ldots, n\}$, such that each element of $X$
occurs at most  once in each
row and at most once in each column.
Hence there are cells in the array that may be empty, but the positions that are
filled have been filled so as to conform with the Latin property of the array.
Let $P$ be a partial Latin square of
order $n$. Then $| P |$ is said to be the {\sf size} of the partial Latin
square and the set of positions ${\cal S}_P=\{(i,j) \mid
(i,j;k)\in P\}$
is said to determine the {\sf shape} of $P$.

A partial Latin square $C$ contained in a Latin square $L$
is said to be {\sf uniquely completable} if $L$ is the
only Latin square of order $n$ with $k$ in position $(i,j)$
for every $(i,j;k) \in C$.
A {\sf critical set} $C$ contained in a Latin square $L$
is a partial Latin square
that is uniquely completable
and no proper subset of $C$ satisfies
this requirement.
The name ``critical set'' and the concept were invented by a statistician,
John Nelder, about 1977, and his ideas were first published 
in a note~{\bf\cite{Nel1}}.  This note posed the problem of giving a formula for the
size of the largest and smallest critical sets for a Latin square of a
given order. Curran and van Rees~{\bf\cite{cv:cs}}, and independently
Smetaniuk~{\bf\cite{smetan}} were the first papers written on the
subject.
See~{\bf\cite{k2}} and~{\bf\cite{vanbate}} for further details.
Let $ \lcs{n}$ denote the size of the {\sf largest critical set}
and $\scs{n}$  the
size of the {\sf smallest critical set} in any Latin square of order $n$.
It was conjectured by Nelder~{\bf\cite{Nel2}}
  that  $ \lcs{n}  =(n^2-n)/2$,
 and by Nelder~{\bf\cite{Nel2}} and also by
 one of the present authors~{\bf\cite{kerman}} and
Bate and van Rees~{\bf\cite{vanbate}}, independently, that
 $\scs{n}  =\lfloor n^2/4 \rfloor$.
The equality for $\lcs{n}$ was shown to be false in 1978, when
Curran and van Rees, {\bf\cite{cv:cs}},
found that $\lcs{4} \geq 7$. Unfortunately, the research over the last
twenty years has not added much
information and in general an upper bound is given by
$n^2-n.$ In this paper we show that
$\lcs{n}  \le n^2-3n+3.$

 In order to validate the construction we require the definition of
 a Latin interchange and an associated  lemma.

 Let $P$ and $P'$ be two
partial Latin squares of the same order and with the same shape.
Then $P$ are $P'$ are said to be {\sf mutually balanced} if the set of
entries in each row (and column) of $P$ are the same as those in the
corresponding row (and column) of $P'$. They are said to be {\sf disjoint}
 if no position in $P'$ contains the same entry as the corresponding position
in $P$.
A {\sf Latin interchange} $I$ is a partial Latin square for
which there exists another partial Latin square $I'$, of the same order,
the same shape and with the property that $I$ and $I'$ are disjoint and
mutually balanced. The partial Latin square $I'$ is said to be a {\sf
disjoint mate} of $I$ (see~{\bf\cite{dha}} and~{\bf\cite{k2}}
for more references).
An example of a Latin interchange and its disjoint mate is given below.
\def\arraystretch{1.0}                 
\begin{center}
$
\begin{array}{|c|c|c|}
\hline
&2& 3  \\ \hline
1&
& 2  \\ \hline
2&3& 1      \\ \hline
\multicolumn{3}{@{\hspace{2pt}}l@{\hspace{2pt}}@{\hspace{0pt}}}{
\hspace*{0.70cm}
I
}
\end{array}
$
\qquad
$
\begin{array}{|c|c|c|}
\hline
&3& 2  \\ \hline
2&
& 1  \\ \hline
1&2& 3      \\ \hline
\multicolumn{3}{@{\hspace{2pt}}l@{\hspace{2pt}}@{\hspace{0pt}}}{
\hspace*{0.70cm}
I'
}
\end{array}
$
\end{center}
The following lemma clarifies the connection between critical sets and Latin
interchanges.
\begin{lemma}
A partial Latin square $C\subseteq L$, of size $s$ and
order $n$, is a critical set for a Latin square $L$ if and only if the
following hold:
\begin{enumerate}
\item[$(1)$]
$C$ contains at least one element of every Latin interchange
that is contained in $L$;
\item[$(2)$]
for each $(i,\, j;\, k)\in C$, there exists a Latin interchange
$I_r$ contained in $L$ such that $I_r\cap C = \{(i,\, j;\, k)\}.$
\end{enumerate}
\end{lemma}
\begin{proof}{
\begin{enumerate}
\item[$(1)$]
If $C$ does not contain an element from some Latin interchange
$I$ in $L$, where $I$ has a disjoint mate $I'$, then $C$ is also a partial
Latin square of $L'= (L\setminus I)\cup I'$. Hence $C$ is not uniquely
completable.
\item[$(2)$]
Since $C$ is a critical set, $C\setminus\{(i,\,j;\,k)\}$
is not uniquely completable. Therefore
$C \setminus \{(i,\,j;\,k)\}$ may be completed in at least two different
ways, thus there exists a Latin interchange $I_r \subseteq L$ such that
$I_r\cap C = \{(i,\, j;\, k)\}.$
\end{enumerate}
\vspace*{-6mm}
}\end{proof}

For a critical set $C$ in a Latin square $L$
we define sets for each row $i$, column $j$ and element $k$. Let \
$R_i = \{ k |\ (i, j; k) \in C \}, \ C_j = \{ k |\ (i, j; k) \in C \}$,
\ and \ $E_k = \{ (i,j) |\ (i, j; k) \in C \}$.
So $R_i$ ($C_j$) is the set of elements which appear in row $i$ (column $j$)
and $E_k$ is the set of positions where the element $k$ appears.
%
\section{The value of $\lcs{n}$  for small $n$}      
In the following table some known values of $\lcs{n}$  are listed for small
values of $n$.  The extra columns are to compare different bounds
discussed in this paper.
$$
\begin{array}{c|cccc}
n & $\lcs{$n$}$  & n^{2}-3n+3 & \lfloor n^{2}-n^{3/2} \rfloor &
\lfloor (1 - (\frac{3}{4})^{log_{2}n}) n^2 \rfloor \\
\hline
 1  &       0  &  1  &  0  &  0  \\
 2  &       1  &  1  &  1  &  1  \\
 3  &       3  &  3  &  3  &  3  \\
 4  &       7  &  7  &  8  &  7  \\
 5  &      11  & 13  & 13  & 12  \\
 6  &      18  & 21  & 21  & 18  \\
 7  & \geq 25  & 31  & 30  & 27  \\
 8  & \geq 37  & 43  & 41  & 37  \\
 9  & \geq 44  & 57  & 54  & 48  \\
10  & \geq 57  & 73  & 68  & 61  \\
\end{array}
$$

The values listed for $\lcs{n}$, 
expect for $n = 5, 7$, 9, and 10,
are given in {\bf\cite{Don1}}.  The value for $n = 5$ and
the bound for $n=7$
were given by A.~Khodkar~{\bf\cite{kho}}.
In the Appendix  we give some examples for the largest
known critical sets for $n = 5, 7$, $9$, and $10$.
The value for $n = 6$ is given in~{\bf\cite{abk}}.

\section{Non-critical sets}    
\setcounter{theorem}{0}
\setcounter{corollary}{0}
\setcounter{lemma}{0}
The following lemma is our main tool in improving the upper
bound on the possible size of $\lcs{n}$.
\begin{lemma}
\label{3.1}
Let $C$ be a critical set for a Latin square $L$ and assume that
there exists $i$
such that $|R_{i}| = n - 1$.
Then the missing element in row $i$ does not occur anywhere in $C$,
and the column corresponding to the missing element is empty.
That is, if ${(i, j; k)} \in L\setminus C$, then
$|C_{j}| =|E_{k}| =  0$.
\end{lemma}
\begin{proof}{
Without loss of generality, let $i=1$ and assume that
$C$ contains the elements $\{ (1,x;x) |\ 1 \leq x \leq n - 1\}$
and that position $(1,n)$ is empty.
Note that the element $n$ may not appear in column $n$ in $C$,
else no element could be placed in position $(1,n)$ of $L$.

By Lemma 1.1 part (2), for each $x$  ($1 \leq x \leq n-1 $)
there exists a Latin interchange $I_x \subseteq L$ such that
$I_x \cap C = \{ (1,x;x)\}$.
Since there is only one empty position in the first row, it follows that 
$\{ (1,x;x), (1,n;n) \} \subseteq I_x$.  Now the interchange $I_x$
has a disjoint mate, say $I'_x$.  In this case since
$ (1,x;n)  \in I'_x$,  for some $r$, $ (r,x; n)  \in I_x$, and since
$|I_x \cap C| = 1$,  $ (r,x;n) \in L\setminus C$.
So $n$ does not occur in column $x$ in $C$.
Since $x$ ranges over all columns from 1 to $n-1$,
$n$ does not occur in $C$ at all.
Therefore $|E_{n}| = 0$.

Also we have $(1,n;x)  \in I'_x$. Thus for some $s$, $ (s,n;x) \in I_x$.
Similarly we have $ (s,n;x)  \notin C$; therefore no element apart from
$n$ may occur in column $n$ in $C$, and we have said that $n$ does
not occur in column $n$ either. Therefore column $n$ is empty.
So $|C_{n}| = 0$.
}\end{proof}
We can generalize Lemma~\ref{3.1} to the following.
\begin{lemma}
\label{3.2}
Let $C$ be a critical set for a Latin square $L$ and assume that
there exists $i$, such that
$|R_i| = n-m$, where
$\{(i, c_1; e_1), (i, c_2; e_2),\dots, (i, c_m; e_m) \}
\subseteq L\setminus C$ \  and
$\{(i,c_{m+1};e_{m+1}), \dots, (i, c_n; e_n) \}\subseteq  C$.
Then we have
\begin{itemize}
\item[$(1)$]
In each of the columns $c_{m+1},c_{m+2}, \dots, c_n$ in $C$, at least
one of the elements $e_1, e_2, \dots, e_m$ is missing.
That is for each $x \in \{c_{m+1}, c_{m+2}, \dots, c_n \}$,
there exists an element $y \in \{e_1, e_2, \dots, e_m \}$, and a row
$r \in \{1,2, 3, \dots, n \}\setminus \{i\}$
such that $(r, x; y)  \in L\setminus C$.
\item[$(2)$]
For each element $e \in \{e_{m+1},e_{m+2}, ..., e_n\}$, we have
a column $c \in \{c_1, c_2, \dots, c_m \}$, from which
this element is missing.
\end{itemize}
\end{lemma}
\begin{proof}{(1)
Without loss of generality we may assume that $i=1$ and $c_j = e_j = j$; for
$j = 1,2, \dots, n$.
For each $x \in \{m+1, m+2, \dots, n \}$,
there exists a Latin interchange $I_x$ such that
$I_x\subseteq L$ and $I_x \cap C = \{ (1, x;x) \}$.
So if $I'_x$ is the disjoint mate of $I_x$ then there exists
$y \in \{1, 2, \dots, m \}$ such that $(1,x;y) \in I'_x$, implying that
there exists $r \in \{ 2, \dots, n \}$ such that $(r,x;y) \in I_x$.
Since $|I_x \cap C| = 1$, $(r,x;y) \in L\setminus C$.

\noindent
(2) Similarly
for each $e \in \{m+1, m+2, \dots, n \}$,
there exists a Latin interchange $I_e$ such that
$I_e\subseteq L$ and $I_e \cap C = \{ (1, e;e) \}$.
So if $I'_e$ is the disjoint mate of $I_e$ then there exists
$c \in \{1, 2, \dots, m \}$ such that $(1,c;e) \in I'_e$, implying that
there exists $s \in \{ 2, \dots, n \}$ such that $(s,c;e) \in I_e$.
Since $|I_e \cap C| = 1$, $(s,c;e) \in L\setminus C$.
}\end{proof}
\begin{theorem}
\label{th3.1}
If $C$ is a uniquely completable partial Latin square of order $n$
completing to the Latin square $L$
with $| C | > n^2-3n+3$, then $C$ is not a critical set.
\end{theorem}
\begin{proof}{
We prove this result by contradiction. Suppose $C$ is a critical set.
Since a critical set in a Latin square of order $n$ can not have $n$
triples whose $i$-th components are the same ($1 \leq i \leq 3$)
(see for example~{\bf\cite{cv:cs}}),
we can assume that any row or column contains at most $n-1$ elements
and any element occurs at most $n-1$ times.

We have three cases to consider.

\noindent
{\bf Case 1 \ }
There exists a row $i$ such that $|R_{i}| = n-1$. Assume that
$(i,j;k) \in L\setminus C.$
Then by Lemma~\ref{3.1}, $|C_{j}| = |E_{k}| =  0$.
Now if there exists $j'$ ($j' \neq j$) such that $|C_{j'}| = n-1$
and $(i', j'; k') \in L\setminus C$,
then we have $|R_{i'}| = 0$. These together imply that
$|C| \le n^2 - (2n-1)-(n-2)=n^2 - 3n + 3$.  Otherwise
$|C_{l}| \leq n-2$, for all $l \neq j$, and $|C_j| = 0$; and thus
$|C|\leq  (n-1)(n -2)  =  n^2 - 3n + 2$.

\noindent
{\bf Case 2 \ }
For all $i$ ($1\le i \le n$) we have, $|R_{i}| \leq n-3$. Then
$|C| \leq n(n - 3) = n^2 - 3n$.

\noindent
{\bf Case 3 \ } For all $i$ ($1\le i \le n$) we have $|R_i| \leq n-2$ and
there exists a row $r$ such that $|R_{r}| = n-2$.
And similarly for all $j$ ($1\le j \le n$) we have $|C_j| \leq n-2$.
Assume that $R_r = \{e_{3},e_{4}, \dots, e_{n}\}$,
and $\{(r, c_{1}; e_{1}), (r, c_{2}; e_{2})\}\subset L\setminus C $.
Then by Lemma~\ref{3.2}
each of the elements $e_{3}, e_4, \dots, e_{n}$
occurs at most once in columns $c_{1}$ and $c_{2}$.   This means
$|C_{c_{1}}| + |C_{c_{2}}| \leq n$. Thus
$|C|\leq  n(n -2) -(n -4) =  n^2 - 3n + 4$.
We will show that $|C| =  n^2 - 3n + 4$ is also impossible. Proof of this fact
is somewhat involved and we need to introduce more notation.

First note that if we consider the conjugate of the Latin square $L$
we may assume that for all $k$ ($1\le k \le n$) we have $|E_k| \leq n-2$.
Let  $f_k = n-2 - |E_k|$. We have $f_k \ge 0$, for all $k$ ($1\le k \le n$).
Assume  $|C| = n^2-3n+4$. Then
\begin{displaymath}
\sum_{k=1}^{n}
f_k = n(n-2)- |C|= n-4.
\end{displaymath}
For each position $(i,j)$, $1 \leq i,j \leq n$, we define $x_{i,j} = |R_i \cup C_j|$.
We have
\begin{displaymath}
\hspace*{-3.3cm}
(*)
\hspace*{3cm}
\sum_{1 \le i,j \le n}
x_{i,j}= n^3 - \sum_{k=1}^{n}(n-|E_k|)^2.
\end{displaymath}
In fact for each position $(i,j)$, $1 \leq i,j \leq n$, we have $x_{i,j} = n$, except when an
element $k$ is missing from {\it both}
row $i$ and column $j$ in $C$.
For each $k$ we have
exactly $(n-|E_k|)^2$  such positions. They are the positions which are in
the $(n-|E_k|) \times (n-|E_k|)$ subsquare obtained from the $n \times n$
array by omitting all the rows and columns containing element $k$ in $C$.
Each such position causes a ``$-1$'' in the summation of the
left hand side of $(*)$.

Note that since $C$ is a critical set, for each  position
$(i,j) \in L\setminus C$, that is for each position in $L$ in which $C$ is
empty, we have \ $x_{i,j} \leq n-1$. Thus
\begin{displaymath}
\begin{array}{ccl}
\frac{1}{|C|}{\displaystyle{\sum_{(i,j)\in C} x_{i,j}}} & = &
\frac{1}{|C|}\Big( (n^3 -{\displaystyle{ \sum_{k=1}^{n}(n-|E_k|)^2)}} -
{\displaystyle{\sum_{(i,j)\in L\setminus C} x_{i,j}}} \Big)    \\
&\ge &
\frac{1}{n^2-3n+4}\Big( (n^3 -{\displaystyle{ \sum_{k=1}^{n}(f_k+2)^2)}} -
(3n-4)(n-1) \Big)  \\
& = &
\frac{1}{n^2-3n+4}{\displaystyle(n^3 - 3n^2-n +12 -{\sum_{k=1}^{n}f^2_k})}. \\
\end{array}
\end{displaymath}
where by $(i,j) \in  C$ we mean a position in $C$  which  is not empty. \\
Since \
${\displaystyle{\sum_{k=1}^{n}f^2_k}} \le
{\displaystyle{(\sum_{k=1}^{n}f_k})^2} = (n-4)^2$, \ thus \
%
$\frac{1}{|C|}{\displaystyle{\sum_{(i,j)\in C} x_{i,j}}} \ge
\frac{n^3 - 3n^2-n +12 -(n-4)^2}{n^2-3n+4} = n-1.$

\noindent
This implies that, either
\begin{enumerate}
\item[(i)]
for some position $(i,j)\in C$ we have $x_{i,j} > n-1$; or
\item[(ii)]
for all $(i,j)\in C$,  $x_{i,j} = n-1$.
\end{enumerate}
The first case is contradictory with $C$ being a critical set.
In the second case if we remove an element $(a,b;e) \in C$, then we have
\begin{itemize}
\item
 $x_{a,b} = n-2$ \ and \
 $x_{a,j} , x_{i,b} \le n-1$, for all \ $(a,j) \ {\rm and} \ (i,b) \in C$; \
 and
 \item
  $x_{i,j} = n-1$; for any other pair $(i,j) \in C$.
\end{itemize}

But if case (ii)
holds, then all of the inequalities that we have above
must be equalities, and
this implies that for every $(i,j) \in L\setminus C$, we have $x_{i,j} = n - 1$.
This follows because we have used the inequality  $x_{i,j}\le  n-1$.
So  $C\setminus\{(a,b;e)\}$ can be completed to $L$, first by completing any
position not in the row $a$ or column $b$, then the positions of row $a$
and column $b$. This is a contradiction.
}\end{proof}
\section{Conjectures and Questions}

There are some conjectures and questions which arise from this research and we
discuss them in this section.

\begin{conjecture}      \ \
$\lcs{n} \leq n^2 - n^{3/2}$.
\end{conjecture}
This is motivated by the proof of Theorem~\ref{th3.1}.
It is analogous to a similar conjecture
made by Brankovic, Horak, Miller, and Rosa, in {\bf\cite{bhmr}}, concerning
the size of the largest premature partial Latin square.
\begin{conjecture} \ \
$\lcs{n} \leq (1 - (\frac{3}{4})^{log_2{n}}) n^{2}$.
\end{conjecture}
This is true for the current known values of $\lcs{n}$. 
It implies that $\lcs{2^{n}} = 4^{n} - 3^{n}$.
This conjecture is based on Stinson and van Rees's result in
{\bf\cite{Sti1}} that $\lcs{2^{n}} \geq 4^{n} - 3^{n}$. We 
postulate that this is an equality.  

\begin{question} \
If $C$ is a critical set of order $n$ and of size $\lcs{n}$,
do there exist $i,j,k$, $1 \leq i,j,k \leq n$, such that
$|R_{i}| = |C_{j}| = |E_{k}| = 0$?  That is, is there
always an empty row, an empty column, and a missing symbol in a critical
set of size $\lcs{n}$?
\end{question}
Evidence for the ``yes'' case in Question 1 is that
every critical set of largest size in Latin
squares of orders 1 to 6 has this property.  Every example in Stinson
and van Rees {\bf\cite{Sti1}} and in Donovan {\bf\cite{Don1}}
where critical sets of largest known size are given, has this property.
All the constructions given for large critical sets given in such articles
as~{\bf\cite{MR99a:05018}},{\bf\cite{ffl}},{\bf\cite{Nel2}} and {\bf\cite{Sti1}}
have this property.
However, the example of a critical set of largest known size in a Latin square
of order 10, given in Appendix 1, does not have this property.

A Latin interchange of size 4 is said to be
an {\sf intercalate}, and the largest number of intercalates in
any Latin square of order $n$ is denoted by $I(n)$ (see {\bf\cite{hw}}).
Below, we ask how $I(n)$, the
maximum number of intercalates in an $n \times n$ Latin square, and
$\lcs{n}$ are related. 

\begin{question} \
If $C$ is a critical set for the Latin square $L$ of order $n$ and size
$\lcs{n}$, does $L$ have $I(n)$ intercalates?
\end{question}

\begin{question} \
If $L$ is a Latin square of order $n$ with $I(n)$ intercalates,
does $L$ contain a critical set $C$ of size $\lcs{n}$?
\end{question}
%

\newpage
\section*{Appendix}
Here we give some examples for the largest
known critical sets for $n = 5, 7$, $9$, and $10$.

\noindent
A critical set of order 5 and size 11:
$$
\begin{latinsq}[1]
\hline  2&&4&3& \\
\hline&&1&2& \\
\hline&2&3&1& \\
\hline  3&1&2&& \\
\hline&&&& \\
\hline
\end{latinsq}
$$
\noindent
A critical set of order 7 and size 25:
$$
\begin{latinsq}[1]
\hline&3&2&1&&5& \\
\hline 6&&3&5&4&1& \\
\hline &6&5&4&3&2& \\
\hline&&4&3&5&& \\
\hline 3&4&1&2&&6& \\
\hline 1&&6&&&3& \\
\hline&&&&&& \\
\hline
\end{latinsq}
$$

In the critical set of order 5, an instance
where $\forall i, |R_{i}| \leq n-2$ has been given
to show that where $C$ is a critical set of size $\lcs{n}$,
it is not necessary to have some  $i,j,k$; $1 \leq i,j,k \leq n$, such
that $|R_{i}| = n - 1$ and  $|C_{j}| = n - 1$ and  $|E_{k}| = n - 1$.

Above, we also gave a similar example for the critical set
of order 7, though it is not known whether $\lcs7 = 25$. And
a critical set of order 9 and size 44 is given below which
also has the same property:
$$
\begin{latinsq}[1]
\hline  1&&  3&&  5&&&  7&\\
\hline&  1&  2&&&&  6&  5& \\
\hline  3&  2&  1&&&  6&  5&  8& \\
\hline&&&  1&&  2&  3&  4& \\
\hline  5&&&  2&  1&  4&  7&  3& \\
\hline&&  5&  3&  2&  1&  4&  6& \\
\hline&  6&  7&  4&  3&&  1&  2& \\
\hline  7&  5&  6&  8&  4&  3&  2&  1& \\
\hline&&&&&&&& \\
\hline
\end{latinsq}
$$

Critical sets of order 9 for all sizes from 20 to 44 inclusive are known
to exist (see~{\bf\cite{Don1}} and {\bf\cite{rwbt}}).

\noindent
A critical set of order 10 and size 57:
$$
\begin{latinsq}
\hline   1   &  & 3   &  & 5   &  & 7   &  & 9 &   \\
\hline    & 1   & 2   &  &  & 5   &  & 6   & 8 &   \\
\hline   3   & 2   & 1   &  &    & 9   & 6   & 7   & 5 &   \\
\hline    &  &  & 1   & 2   & 3   &  & 8   & 4 &   \\
\hline   5   &  &  & 2   & 1   &  10  & 4   & 3   &&   \\
\hline    & 5   & 9   & 3   &  10  & 1   & 2   &  & 6 &   \\
\hline   7   &  & 6   &  & 4   & 2   & 1   & 5   & 3 &   \\
\hline    & 6   & 7   & 8   & 3   &  & 5   & 1   & 2 &   \\
\hline   9   & 8   & 5   & 4   &    & 6   & 3   & 2   & 1 &   \\
\hline    &  &  &  &    &  &  &  &&   \\
\hline
\end{latinsq}
$$

Critical sets of order 10 for all sizes from 25 to 57 inclusive are known
to exist (see {\bf\cite{Don1}}, {\bf\cite{closing}}, and {\bf\cite{rwbt}}).
\section*{Acknowledgements}
We thank Diane Donovan for her very useful comments on this paper.
We also appreciate M. Mahdian, P.J. Owens, and G.H.J. van Rees
for reading the preprint and their comments.
Richard Bean wishes to thank the Australian Government
for the support of an Australian Postgraduate Award and an Australian
Research Council Large Grant  A49937047.
E.S. Mahmoodian appreciates the hospitality
of the Department of Mathematics at the University of
Queensland, while working on this paper, his research is supported by
Australian Research Council Grant A69701550.


\end{document}